# THE MINIMAL ENTROPY MARTINGALE MEASURE FOR GENERAL BARNDORFF-NIELSEN/SHEPHARD MODELS[1]

### By Thorsten Rheinländer and Gallus Steiger

*London School of Economics and ETH Zürich*

We determine the minimal entropy martingale measure for a general class of stochastic volatility models where both price process and volatility process contain jump terms which are correlated. This generalizes previous studies which have treated either the geometric Lévy case or continuous price processes with an orthogonal volatility process. We proceed by linking the entropy measure to a certain semilinear integro-PDE for which we prove the existence of a classical solution.

**1. Introduction.** The main contribution of this paper is the calculation of the minimal entropy martingale measure (MEMM) for a general class of stochastic volatility models as explicitly as possible in terms of the parameters of the market model. This idea of explicit description of optimal martingale measures started with the minimal martingale measure of Föllmer and Schweizer [12], followed by the minimal Hellinger martingale measure of Grandits [15] and the minimal entropy-Hellinger martingale measure of Choulli and Stricker [8].

Our study of the MEMM encompasses the simpler cases where either the dynamics of the risky asset is modeled as a geometric Lévy process or the price process is continuous with an orthogonal pure jump volatility process. These cases, as will be discussed below, have been studied separately and with different methods. Our approach presents a unifying framework which moreover covers models like the Barndorff-Nielsen and Shephard (BN–S) model where both price process and volatility process contain jump terms which are correlated. It turns out that, due to the correlation, this general case is much more difficult and can be considered a nontrivial mixture of the two cases studied previously.

Received March 2005; revised November 2005.

[1]Supported by NCCR Financial Valuation and Risk Management.

*AMS 2000 subject classifications.* 28D20, 60G48, 60H05, 91B28.

*Key words and phrases.* Relative entropy, martingale measures, stochastic volatility.

Asset process models driven by nonnormal Lévy processes date back to the work of Mandelbrot [22]. More recently, rather complex models like the stochastic volatility model of Barndorff-Nielsen and Shephard [1] have been developed. This model is constructed via a jump-diffusion price process together with a mean reverting, stationary volatility process of Ornstein–Uhlenbeck type driven by a subordinator (i.e., an increasing Lévy process). Moreover, the negative correlation between price process and volatility process in this model allows us to deal with the so-called leverage problem, that is, for equities, a fall in price level is typically associated with an increase in volatility.

One main reason for the use of Lévy-driven asset models is the flexibility they allow when fitting a model to observed asset prices. However, the corresponding financial market is then typically incomplete, resulting in the existence of multiple equivalent martingale measures. A standard approach is to identify an optimal martingale measure on the basis of the utility function of the investor; see [20]. In this paper, we consider the exponential utility function which corresponds via an asymptotic utility indifference approach to taking the MEMM as pricing measure [3, 4, 9].

In case the price process is an exponential Lévy process, the MEMM has been calculated by several authors in varying degrees of generality (e.g., [7, 11, 14, 23]). Grandits and Rheinländer [16] and Benth and Meyer-Brandis [6] determine the MEMM in stochastic volatility models where the price process is driven by a Brownian motion $B$, and the volatility process may contain jump terms and is orthogonal to $B$. Still assuming a continuous price process, Becherer [3] considers a model with interacting Itô and point processes.

With respect to the BN–S model with leverage effect, Nicolato and Venardos [24] analyze the class of all equivalent martingale measures, with a focus on the subclass of structure-preserving martingale measures (i.e., the price process is also of BN–S-type under those martingale measures). In the case of exponential Lévy processes, the asset process under the MEMM is again an exponential Lévy process (see in particular [11]), but one major implication of the results in this paper is that the volatility process in the BN–S model in general no longer has independent increments under the MEMM. Therefore, only considering the class of structure-preserving martingale measures seems to be too narrow an approach, especially in the context of exponential utility maximization.

The paper is structured as follows. Section 2 introduces our setup and the martingale approach for determining the MEMM in case of a general Lévy process-driven asset model. In Section 3 we consider a general class of stochastic volatility models. We derive the structure of the MEMM by linking it to the solution of a certain semi-linear integro-PDE, a unique classical solution of which is shown to exist. We conclude this paper in Section 4 by applying this result to the two extreme cases—(1) price process is a Lévy process and (2) price process is continuous with an orthogonal stochastic volatility process—as well as to the BN–S model. The latter case presents an additional technical difficulty since the volatility process is unbounded. This issue has been resolved in [27].

The present approach has been influenced by the martingale duality approach in Rheinländer [26] where the MEMM was linked to the solution of a certain equation in the case of a filtration where all martingales are continuous. This has been applied in [18] and [26] to stochastic volatility models driven by Brownian motions. The presence of jumps, however, calls for more general techniques. Our method was inspired by Becherer's [2] approach which considers interacting systems of semi-linear PDEs.

## 2. Preliminaries and general results.
We start with some general assumptions which hold throughout the paper. Let $(\Omega, \mathcal{F}, \mathbb{F}, P)$ be a filtered probability space and $T$ some fixed finite time horizon. We assume that $\mathcal{F}_0$ is trivial and that $\mathcal{F} = \mathcal{F}_T$. The filtration $\mathbb{F} = (\mathcal{F}_t)_{0 \le t \le T}$ fulfills the usual conditions and is generated by a Lévy process $Y$ where $Y^c$ $(Y^d)$, $\mu_Y$, and $\nu_Y(dx, dt) = \nu(dx) \, dt$ denote its continuous (discontinuous) martingale part, the jump measure and its compensator, respectively. For simplicity, we assume that $\langle Y^c \rangle_t = t$. We refer to [19] with respect to the notation used in this paper. In particular, $G_{\mathrm{loc}}(\mu_Y)$ is defined in [19], Definition II.1.27.

REMARK 2.1. By Jacod and Shiryaev [19], Theorem III.4.34, we have the following *representation property*: every $(P, \mathbb{F})$-local martingale $M$ can be written as

$$M = M_0 + \int H \, dY^c + W * (\mu_Y - \nu_Y)$$

for some $H \in L^2_{\mathrm{loc}}(Y^c)$, $W \in G_{\mathrm{loc}}(\mu_Y)$.

We denote by $S$ an $\mathbb{F}$-adapted, locally bounded semimartingale (modeling the price process of a risky asset), which has the following canonical decomposition:

$$S = S_0 + M + A,$$

where $M$ is a locally bounded local martingale with $M_0 = 0$ and $A$ is a process of locally finite variation. By the representation property, we write $M$ as

$$M = M^c + M^d = \int \sigma^M \, dY^c + W^M(x) * (\mu_Y - \nu_Y),$$

where $M^c$ and $M^d$ are the continuous and the discontinuous parts of the local martingale $M$, respectively, $\sigma^M$ is predictable and $W^M \in G_{\mathrm{loc}}(\mu_Y)$. Moreover, we assume that the asset price process $S$ satisfies the following:

ASSUMPTION 2.2 (Structure condition). There exists a predictable process $\lambda$ satisfying

$$A = \int \lambda \, d\langle M \rangle,$$

with

$$K_T := \int_0^T \lambda_s^2 \, d\langle M \rangle_s < \infty, \qquad P\text{-a.s.}$$

DEFINITION 2.3. Let $\mathcal{V}$ be the linear subspace of $L^\infty(\Omega, \mathcal{F}, P)$ spanned by the elementary stochastic integrals of the form $f = h(S_{T_2} - S_{T_1})$, where $0 \leq T_1 \leq T_2 \leq T$ are stopping times such that the stopped process $S^{T_2}$ is bounded and $h$ is a bounded $\mathcal{F}_{T_1}$-measurable random variable. A *martingale measure* is a probability measure $Q \ll P$ with $E[\frac{dQ}{dP} f] = 0$ for all $f \in \mathcal{V}$.

We denote by $\mathcal{M}$ the set of all martingale measures for $S$ and by $\mathcal{M}^e$ the subset of $\mathcal{M}$ consisting of probability measures which are equivalent to $P$. Here and in the sequel, we identify measures with their densities. Note that, as $S$ is locally bounded, a probability measure $Q$ absolutely continuous to $P$ is in $\mathcal{M}$ if and only if $S$ is a local $Q$-martingale.

DEFINITION 2.4. The relative entropy $I(Q, R)$ of the probability measure $Q$ with respect to the probability measure $R$ is defined as

$$I(Q, R) = \begin{cases} E_R\left[\dfrac{dQ}{dR} \log \dfrac{dQ}{dR}\right], & \text{if } Q \ll R, \\ +\infty, & \text{otherwise.} \end{cases}$$

It is well known that $I(Q, R) \geq 0$ and that $I(Q, R) = 0$ if and only if $Q = R$.

DEFINITION 2.5. The minimal entropy martingale measure $Q^E$, also abbreviated MEMM in what follows, is the solution of

$$\min_{Q \in \mathcal{M}} I(Q, P).$$

Theorems 1, 2 and Remark 1 of [13], as well as the fact that $\mathcal{V} \subset L^\infty(P)$, yield the following:

THEOREM 2.6 ([13]). *If there exists $Q \in \mathcal{M}^e$ such that $I(Q, P) < \infty$, then the minimal entropy martingale measure exists, is unique and moreover is equivalent to $P$.*

Let us restate the following criterion for a martingale measure to coincide with the MEMM:

THEOREM 2.7 ([16]). *Assume there exists a $Q \in \mathcal{M}^e$ with $I(Q, P) < \infty$. Then $Q^*$ is the minimal entropy martingale measure if and only if there exists a constant $c$ and an $S$-integrable predictable process $\phi$*

$$\frac{dQ^*}{dP} = \exp\left(c + \int_0^T \phi_t \, dS_t\right) \tag{2.1}$$

*such that $E_Q[\int_0^T \phi_t \, dS_t] = 0$ for all $Q \in \mathcal{M}^e$ with finite relative entropy.*

REMARK 2.8. Based on the above results, we will pursue the following strategy to determine the MEMM. We first find some candidate measure $Q^*$ which can be represented as in (2.1). To verify that $Q^*$ is indeed the entropy minimizer, we then proceed in three steps, showing that:

1. $Q^*$ is an equivalent martingale measure;
2. $I(Q^*, P) < \infty$;
3. $\int \phi \, dS$ is a true $Q$-martingale for all $Q \in \mathcal{M}^e$ with finite relative entropy.

This martingale approach yields a necessary equation for $\phi$ and $c$:

THEOREM 2.9. *Assume that the MEMM $Q^*$ exists. The strategy $\phi$ and the constant $c$ in (2.1) satisfy the equation*

$$
\begin{aligned}
c + \int_0^T &\left[\tfrac{1}{2}(\sigma_t^L - \lambda_t \sigma_t^M)^2 + \phi_t \lambda_t (\sigma_t^M)^2 + \phi_t \lambda_t \int_{\mathbb{R}} (W_t^M(x))^2 \nu(dx)\right] dt \\
&= \int_0^T (\sigma_t^L - (\phi_t + \lambda_t)\sigma_t^M) \, dY_t^c \\
&\quad + ((W^L(x) - (\phi + \lambda)W^M(x)) * (\mu_Y - \nu_Y))_T \\
&\quad + ((\log(1 - \lambda W^M(x) + W^L(x)) + \lambda W^M(x) - W^L(x)) * \mu_Y)_T
\end{aligned}
\tag{2.2}
$$

*with predictable processes $\sigma^L \in L_{\mathrm{loc}}^2(Y^c)$ and $W^L \in G_{\mathrm{loc}}(\mu_Y)$ such that*

$$\sigma_t^M \sigma_t^L + \int_{\mathbb{R}} W_t^M(x) W_t^L(x) \nu(dx) = 0 \qquad \forall t \in [0, T]. \tag{2.3}$$

PROOF. By Girsanov's theorem together with the structure condition, the density process $Z = (Z_t)$ of $Q^*$ is a stochastic exponential of the form

$$Z = \mathcal{E}\left(-\int \lambda \, dM + L\right),$$

where $L$ and $[M, L]$ are local $P$-martingales. Using the representation property, let us write the local martingale $L$ in the following way:

$$L = \int \sigma^L \, dY^c + W^L(x) * (\mu_Y - \nu_Y), \tag{2.4}$$

for some $\sigma^L \in L^2_{\mathrm{loc}}(Y^c)$, $W^L \in G_{\mathrm{loc}}(\mu_Y)$. We therefore get

$$[M, L] = \int_0^{\cdot} \sigma_s^M \sigma_s^L \, ds + W^M(x) W^L(x) * \mu_Y.$$

Furthermore, we know from [10], VII.39, that the predictable bracket process

$$\langle M, L \rangle = \int_0^{\cdot} \sigma_s^M \sigma_s^L \, ds + W^M(x) W^L(x) * \nu_Y$$

exists, since $M$ is locally bounded. However, $\langle M, L \rangle$ is equal to zero since $[M, L]$ is a local martingale. Therefore, condition (2.3) holds. We now apply Itô's formula to $\log Z$ to get, for $t \in [0, T]$, that

$$
\begin{aligned}
\log Z_t &= \int_0^t \frac{1}{Z_{s-}} \, dZ_s - \frac{1}{2} \int_0^t \frac{1}{Z_{s-}^2} \, d\langle Z^c \rangle_s \\
&\quad + \sum_{s \le t} \left( \log Z_s - \log Z_{s-} - \frac{1}{Z_{s-}} \Delta Z_s \right) \\
&= - \int_0^t \lambda_s \, dM_s + L_t - \frac{1}{2} \int_0^t \lambda_s^2 \, d\langle M^c \rangle_s + \int_0^t \lambda_s \, d\langle M^c, L^c \rangle_s - \frac{1}{2} \langle L^c \rangle_t \\
&\quad + \sum_{s \le t} \left( \log \frac{Z_s}{Z_{s-}} + \Delta \int_0^s \lambda \, dM - \Delta L_s \right) \\
&= \int_0^t (\sigma_s^L - \lambda_s \sigma_s^M) \, dY_s^c - \frac{1}{2} \int_0^t (\lambda_s \sigma_s^M - \sigma_s^L)^2 \, ds \\
&\quad + ((W^L(x) - \lambda W^M(x)) * (\mu_Y - \nu_Y))_t \\
&\quad + ((\log(1 - \lambda W^M(x) + W^L(x)) + \lambda W^M(x) - W^L(x)) * \mu_Y)_t.
\end{aligned}
$$

Moreover, due to Theorem 2.7, at the time horizon we have

$$
\begin{aligned}
\log Z_T &= c + \int_0^T \phi_t \, dS_t \\
&= c + \int_0^T \phi_t \sigma_t^M \, dY_t^c + (\phi W^M(x) * (\mu_Y - \nu_Y))_T \\
&\quad + \int_0^T \left( \phi_t \lambda_t (\sigma_t^M)^2 + \phi_t \lambda_t \int_{\mathbb{R}} (W_t^M(x))^2 \nu(dx) \right) dt.
\end{aligned}
$$

We arrive at equation (2.2) by combining the two equations above. $\quad \square$

COROLLARY 2.10. *Equation* (2.2) *in Theorem* 2.9 *is fulfilled once the following conditions are satisfied:*

(i) $|W^L(x) - (\phi + \lambda) W^M(x)| * \mu_Y \in \mathcal{A}_{\mathrm{loc}}^+$.

(ii) It holds that

$$
\begin{aligned}
c + \int_0^T & [\tfrac{1}{2}(\sigma_t^L - \lambda_t \sigma_t^M)^2 + \phi_t \lambda_t (\sigma_t^M)^2]\, dt \\
& + \int_0^T \int_{\mathbb{R}} (W_t^L(x) - (\phi_t + \lambda_t) W_t^M(x) + \phi_t \lambda_t (W_t^M(x))^2)\, \nu(dx)\, dt \\
& = \int_0^T (\sigma_t^L - (\phi_t + \lambda_t)\sigma_t^M)\, dY_t^c \\
& \qquad + ((\log(1 - \lambda W^M(x) + W^L(x)) - \phi W^M(x)) * \mu_Y)_T.
\end{aligned}
\tag{2.5}
$$

PROOF. By Jacod and Shiryaev [19], Proposition II.1.28, condition (i) implies that we can write

$$
\begin{aligned}
(W^L(x) - (\phi + \lambda) W^M(x)) * (\mu_Y - \nu_Y) = {} & (W^L(x) - (\phi + \lambda) W^M(x)) * \mu_Y \\
& - (W^L(x) - (\phi + \lambda) W^M(x)) * \nu_Y.
\end{aligned}
$$

Taking this into account, equation (2.2) reduces to the simpler equation (2.5). $\square$

Once we have, by solving (2.2) and (2.3), found a candidate martingale measure, we still have to carry out the verification procedure outlined above. We will need the following lemma, which is a generalization of the Novikov condition to discontinuous processes:

LEMMA 2.11 ([21]). *Let $N$ be a locally bounded local $P$-martingale. Let $Q$ be a measure defined by*

$$
\frac{dQ}{dP}\Big|_{\mathcal{F}_t} = Z_t = \mathcal{E}(N)_t,
$$

*where $\Delta N > -1$. If the process*

$$
U_t = \tfrac{1}{2}\langle N^c \rangle_t + \sum_{s \le t} \{(1 + \Delta N_s)\log(1 + \Delta N_s) - \Delta N_s\}
\tag{2.6}
$$

*belongs to $\mathcal{A}_{\mathrm{loc}}$, and therefore has a predictable compensator $B_t$ and, in addition,*

$$
E[\exp B_T] < \infty,
\tag{2.7}
$$

*then $Q$ is an equivalent probability measure.*

Finally, to cope with item 3 of our approach described in Remark 2.8, we mention the following result:

LEMMA 2.12 ([26]). *Let $Q$ be an equivalent martingale measure with finite relative entropy, and let $\int \psi \, dS$ be a local $Q$-martingale. Then $\int \psi \, dS$ is a true $Q$-martingale if, for some $\beta > 0$ small enough, $\exp\{\beta \int_0^T \psi_t^2 \, d[S]_t\}$ is $P$-integrable.*

## 3. A general jump-diffusion model.

Let us consider a class of stochastic volatility models with asset prices of the following type:

$$(3.1) \quad \frac{dS_t}{S_{t-}} = \eta^M(t, V_t) \, dt + \sigma^M(t, V_t) \, dY_t^c + d(W^M(\,\cdot\,, V_-, x) * (\mu_Y - \nu_Y))_t,$$

$$(3.2) \quad dV_t = \eta^V(t, V_t) \, dt + d(W^V(\,\cdot\,, V_-, x) * \mu_Y)_t,$$

where $V$ is defined on some interval $E \subset \mathbb{R}$. In the notation we will often suppress the dependence on $V$ of the various processes. Our basic assumptions are as follows:

ASSUMPTION 3.1.

1. The coefficient $\eta^V$ is differentiable in $y$ (corresponding to the "V-coordinate") with bounded continuous partial derivative and is locally Lipschitz-continuous in $t$. $W^V$ is differentiable in $y$ with bounded derivative and continuous in $t$.

2. The coefficients $\eta^M$, $\sigma^M$ and $W^M$ are locally Lipschitz-continuous in $t$ and differentiable in $y$ with bounded derivative. Furthermore, $\eta^M$ is positive, $\sigma^M$ is positive and uniformly bounded away from zero on $[0, T] \times E$ and $W_t^M(y, \cdot) : \text{supp}(\nu) \to (-1, \infty) \in l^\infty(\text{supp}(\nu)) \cap l^1(\text{supp}(\nu))$, uniformly in $t$.

3. The functions $W^M$ and $W^V$ are in $G_{\text{loc}}(\mu)$.

4.

$$\widehat{\lambda} := \frac{\eta^M}{(\sigma^M)^2 + \int (W^M(x))^2 \nu(dx)}$$

   is uniformly bounded on $[0, T] \times E$.

5. We have $\int |W^V(x)| \nu(dx) < \infty$.

REMARK 3.2. By Protter [25], Theorem V.38 and the remark following it, Assumptions 3.1.1–3.1.3 ensure that there exists a unique solution $(S, V)$ to equations (3.1) and (3.2) which does not explode in $[0, T]$.

Let us turn to our basic equation (2.5). The functions $\sigma^M$ and $W^M(x)$ of Section 2 correspond now, with a slight abuse of notation, to $\sigma^M S_-$ and $W^M(x) S_-$, respectively. Moreover, we set $\widehat{\lambda} := \lambda S_-$ and $\widehat{\phi} := \phi S_-$. We denote

$$\Delta u_t = \Delta u_t(y, x) := u(t, y + W^V(t, y, x)) - u(t, y),$$

and work with the ansatz that there exists a sufficiently smooth function $u$ such that

$$(3.3) \qquad (\log(1 - \widehat{\lambda} W^M + W^L) - \widehat{\phi} W^M)(t, V_{t-}, x) = \Delta u_t(V_{t-}, x),$$

that is, the jumps of the right-hand side of (2.5) correspond to the jumps of some function $u$ along the paths of process $V$. In addition, we set

$$(3.4) \qquad u(T, \cdot) = 0 \qquad \text{on } E.$$

With this ansatz we can write, using Itô's formula,

$$([\log(1 - \widehat{\lambda} W^M(x) + W^L(x)) - \widehat{\phi} W^M(x)] * \mu_Y)_T$$
$$= \sum_{0 < t \le T} \{u(t, V_t) - u(t, V_{t-})\}$$
$$= -u(0, V_0) - \int_0^T \left( \frac{\partial}{\partial t} u(t, V_t) + \eta_t^V \frac{\partial}{\partial V} u(t, V_t) \right) dt.$$

We may therefore rewrite equation (2.5) as

$$c + u(0, V_0)$$
$$= -\int_0^T \left[ \frac{1}{2}(\sigma_t^L - \widehat{\lambda}_t \sigma_t^M)^2 + \widehat{\phi}_t \widehat{\lambda}_t (\sigma_t^M)^2 \right.$$
$$(3.5) \qquad\qquad + \frac{\partial}{\partial t} u(t, V_t) + \eta_t^V \frac{\partial}{\partial V} u(t, V_t)$$
$$\left. + \int (W_t^L(x) - (\widehat{\phi}_t + \widehat{\lambda}_t) W_t^M(x) + \widehat{\phi}_t \widehat{\lambda}_t (W_t^M(x))^2) \nu(dx) \right] dt$$
$$+ \int_0^T [\sigma_t^L - (\widehat{\phi}_t + \widehat{\lambda}_t) \sigma_t^M] \, dY_t^c.$$

A solution to this problem might be to require that

$$(3.6) \qquad \frac{1}{2}(\sigma^L - \widehat{\lambda} \sigma^M)^2 + \widehat{\phi} \widehat{\lambda}(\sigma^M)^2 + \frac{\partial}{\partial t} u(\cdot, V) + \eta^V \frac{\partial}{\partial V} u(\cdot, V)$$
$$+ \int (W^L(x) - (\widehat{\phi} + \widehat{\lambda}) W^M(x) + \widehat{\phi} \widehat{\lambda}(W^M(x))^2) \nu(dx) = 0,$$

together with (3.4) and

$$(3.7) \qquad c = -u(0, V_0), \qquad \sigma^L = (\widehat{\phi} + \widehat{\lambda}) \sigma^M.$$

Let us introduce $u_t := u(t, \cdot) \colon E \to \mathbb{R}$ and

$$g^y(t, u_t) := \frac{1}{2}(\sigma_t^L(y) - \widehat{\lambda}_t(y) \sigma_t^M(y))^2 + \widehat{\phi}_t(y) \widehat{\lambda}_t(y)(\sigma_t^M(y))^2$$
$$(3.8) \qquad + \int (W_t^L(y, x) - (\widehat{\phi}_t(y) + \widehat{\lambda}_t(y)) W_t^M(y, x)$$
$$+ \widehat{\phi}_t(y) \widehat{\lambda}_t(y)(W_t^M(y, x))^2) \nu(dx).$$

Provided that $\widehat{\phi}_t$, $\sigma_t^L$ and $W_t^L(x)$ are functions of $u_t$, (3.6) is an integro-PDE for $u$ of the form

$$(3.9) \qquad \frac{\partial}{\partial t} u(t,y) + \eta_t^V \frac{\partial}{\partial y} u(t,y) + g^y(t, u_t) = 0,$$

$$(3.10) \qquad\qquad\qquad u(T,y) = 0 \qquad \text{for all } y \in E.$$

By equation (3.7) together with condition (2.3), we get

$$(3.11) \qquad \widehat{\phi} = -\widehat{\lambda} - \frac{\int W^M(x) W^L(x) \nu(dx)}{(\sigma^M)^2},$$

which, by equation (3.3), leads to (suppressing the $t$ and $y$ variables)

$$(3.12) \quad \begin{aligned} &\exp\left\{ \Delta u(x) - \left[ \widehat{\lambda} + \frac{\int W^M(z) W^L(z) \nu(dz)}{(\sigma^M)^2} \right] W^M(x) \right\} \\ &= 1 - \widehat{\lambda} W^M(x) + W^L(x). \end{aligned}$$

To make this intuitive approach rigorous, we shall proceed as follows. We show in Corollary 3.4 below that each $u \in \mathcal{C}_b([0,T] \times E)$ gives via $\Delta u$ a unique bounded function $W^L$ solving (3.12). We then define $\widehat{\phi}$ as in (3.11), $\sigma^L$ as in (3.7) and $g^y$ as in (3.8). In Theorem 3.8 below it is then shown that there exists a classical solution to the integro-PDE (3.9), (3.10). Finally, we provide the verification results in Theorem 3.9.

For the discussion of equation (3.12) we first provide a preparatory result:

**Lemma 3.3.** *Let $\beta > 0$, $f \in l^\infty(\mathrm{supp}(\nu)) \cap l^1(\mathrm{supp}(\nu))$, the set of bounded and integrable functions from $\mathrm{supp}(\nu)$ into $\mathbb{R}$, and $k$ be a function on $\mathrm{supp}(\nu)$ which is bounded from above. Then, the function $\varphi \colon \mathrm{supp}(\nu) \to \mathbb{R}$, given as*

$$\varphi(x) = \exp\left\{ k(x) - \beta f(x) \int f(z) \varphi(z) \nu(dz) \right\},$$

*is well defined and bounded.*

For the proof see the Appendix.

**Corollary 3.4.** *Let Assumption 3.1 hold and let $u$ be defined on $([0, T] \times E)$ such that $\Delta u$ is uniformly bounded from above. Then, $u$ uniquely defines a function*

$$W^L = W^L(u) \colon [0,T] \times E \times \mathbb{R} \to \mathbb{R}$$

*which fulfills equation (3.12). $W^L$ and, therefore, also $\widehat{\phi}$ and $\sigma^L$, is uniformly bounded in $(t,y) \in [0,T] \times E$ and, moreover, $W_t^L(y, \cdot) \in l^\infty(\mathrm{supp}(\nu)) \cap l^1(\mathrm{supp}(\nu))$ for all $(t,y) \in [0,T] \times E$.*

Proof.    Introducing

$$(3.13) \qquad \varphi(x) := (W_t^L(u_t))(x) - \widehat{\lambda} W_t^M(x) + 1,$$

we may write equation (3.12) pointwise in $t \in [0, T]$ in the form

$$\varphi(x) = \exp\left\{k(x) - \beta f(x) \int f(z)\varphi(z)\nu(dz)\right\}$$

with

$$f(x) := W_t^M(x),$$

$$k(x) := \Delta u_t(x) - W_t^M(x)\left[\widehat{\lambda}_t\left(1 + \frac{\int (W_t^M(z))^2\nu(dz)}{(\sigma_t^M)^2}\right) - \frac{\int W_t^M(z)\nu(dz)}{(\sigma_t^M)^2}\right],$$

$$\beta := \frac{1}{(\sigma_t^M)^2}.$$

Since $\Delta u_t$ is bounded from above, we have that $k$ is bounded from above, by Assumption 3.1, and we may apply Lemma 3.3. By the definition of $\varphi$ in (3.13), it follows directly that $W^L$ fulfills equation (3.12). $W^L$ is even uniformly bounded in $(t, y)$ since $\Delta u$ is uniformly bounded from above. Finally, we get $W_t^L(y, \cdot) \in l^1(\text{supp}(\nu))$ from a Taylor expansion together with our assumption that $\int |W^V(x)|\nu(dx)$.    □

The function $W_t^L(y, \cdot)$, seen as a function of $u_t$,

$$\mathcal{C}_b(E) \to l^1(\text{supp}(\nu)),$$

$$u_t \mapsto (W_t^L(u_t))(y, \cdot),$$

is not uniformly Lipschitz-continuous. However, we can ensure this property by restricting the space $\mathcal{C}_b(E)$ to the set

$$\mathcal{C}_b^Q(E) := \{v \in \mathcal{C}_b(E), \|v\|_\infty \leq Q\},$$

with a constant $Q > 0$. In fact, we obtain the following:

Lemma 3.5.    *For $(t, y) \in [0, T] \times E$ fixed,*

$$W_t^L(y, \cdot) : \mathcal{C}_b^Q(E) \to l^1(\text{supp}(\nu))$$

*is Lipschitz-continuous, uniformly with respect to $t \in [0, T]$, and with a Lipschitz constant independent of $y$.*

For the proof see the Appendix.

We turn now to the existence of a solution for the integro-PDE (3.9)–(3.10). The following two theorems provide some general existence results:

THEOREM 3.6. *Let $E \subset \mathbb{R}$ be some interval. For $(t, z) \in [0, T] \times E$, consider*

$$(3.14) \qquad Z_{\cdot}^{t,z} = z + \int_t^{\cdot} b(u, Z_u^{t,z}) \, du,$$

*for a continuous process $b \colon [0, T] \times E \to \mathbb{R}$, such that $Z^{t,z}$ stays in $E$.*

*Let us consider the partial differential equation with boundary condition:*

$$(3.15) \qquad \frac{\partial}{\partial t} u(t, z) + b(t, z) \frac{\partial}{\partial z} u(t, z) + g^z(t, u_t) = 0,$$

$$(3.16) \qquad\qquad u(T, z) = h(z) \qquad \forall z \in E,$$

*for which we shall assume:*

(a-1) *$b$ is locally Lipschitz-continuous.*

(a-2) *$g \colon [0, T] \times \mathcal{C}_b(E) \to \mathcal{C}_b(E)$ is a Lipschitz-continuous function in $v \in \mathcal{C}_b(E)$, uniformly in $t$. That is, there exists a constant $L < \infty$ such that*

$$\|g(t, v_1) - g(t, v_2)\|_\infty \leq L \|v_1 - v_2\|_\infty \qquad \forall t \in [0, T], v_1, v_2 \in \mathcal{C}_b(E).$$

(a-3) *$h \colon E \to \mathbb{R} \in \mathcal{C}_b(E)$.*

*Then, there exists a unique solution $\widehat{u} \in C_b([0, T] \times E)$ which solves the boundary problem* (3.15)–(3.16) *in the sense of distributions. It can be written as*

$$\widehat{u}(t, z) = h(Z_T^{t,z}) + \int_t^T g^{Z_s^{t,z}}(s, \widehat{u}_s) \, ds.$$

For the proof see the Appendix.

Existence of a strong solution can be ensured in the following special case:

THEOREM 3.7. *Let us assume that all conditions of Theorem 3.6 are fulfilled. Let us further assume that $E \subset \mathbb{R}$ is compact and that the following hold true:*

(b-1) *$b$ has a uniformly bounded, continuous derivative $\frac{\partial}{\partial z} b$.*

(b-2) *For any $v \in C_1^b(E)$, $g^z(t, v)$ is differentiable in $z$ with $\frac{\partial}{\partial z} g^z(t, v) = \widehat{g}^z(t, \frac{\partial}{\partial z} v)$ for some suitable continuous function $\widehat{g}$, such that:*

• *there exist some constants $L, K$ such that we may write*

$$(3.17) \qquad \|\widehat{g}(s, v_s)\|_\infty \leq L \|v_s\|_\infty + K,$$

• *for any $R > 0$, $\widehat{g}$ is uniformly continuous on $[0, T] \times M \times E$ with $M = C_b^R(E)$,*

(b-3) $h \in \mathcal{C}_b^1(E)$.

*Then, the weak solution $\hat{u} \in \mathcal{C}_b([0,T] \times E)$ is differentiable in the space variable and, therefore, it is also the strong solution to the boundary problem* (3.15)–(3.16).

For the proof see the Appendix.

Let us apply this result to $g^y(t, u_t)$ having the form (3.8). In this case, $g^y(t, u_t)$ does not have to be Lipschitz-continuous. However, using a truncation argument we get the following result:

THEOREM 3.8. *Let Assumption* 3.1 *hold and let $g^y(t, u_t)$ be of the form* (3.8). *Let $E$ be a compact interval such that $\sigma^M$ is uniformly bounded on $[0,T] \times E$. Then there is a classical solution $\hat{u} \in \mathcal{C}_b^{1,1}([0,T] \times E)$ to the integro-PDE*

$$(3.18) \qquad \frac{\partial}{\partial t} u(t,y) + \eta_t^V \frac{\partial}{\partial y} u(t,y) + g^y(t, u_t) = 0$$

*with boundary condition*

$$(3.19) \qquad u(T, y) = 0.$$

$\hat{u}$ *satisfies*

$$(3.20) \qquad \hat{u}(t,y) = \int_t^T g^{\hat{V}_s^{t,y}}(s, \hat{u}_s)\, ds$$

*with*

$$(3.21) \qquad d\hat{V}_s^{t,y} = \eta^V(s, \hat{V}_s^{t,y})\, ds$$

*and $\hat{V}_t^{t,y} = y$.*

PROOF. Let us rewrite (3.8) using (3.7) and (3.11) as

$$
\begin{aligned}
g(\cdot, v) = \frac{1}{2} \Bigg[ & \left( \frac{\int W^M(x) W^L(x) \nu(dx)}{\sigma^M} \right)^2 - \hat{\lambda}^2 (\sigma^M)^2 \Bigg] \\
& + \frac{\int W^M(x) \nu(dx) - \hat{\lambda} \int (W^M(x))^2 \nu(dx)}{(\sigma^M)^2} \\
& \times \int W^M(x) W^L(x) \nu(dx) \\
& + \int W^L(x) \nu(dx) - \hat{\lambda}^2 \int (W^M(x))^2 \nu(dx),
\end{aligned}
$$

(3.22)

which is in general not Lipschitz-continuous. We circumvent this problem by introducing a truncating, auxiliary function $\tilde{g}$. We will show that the weak solution $\widehat{u} \in \mathcal{C}_b([0,T] \times \mathbb{R})$ to the integro-PDE

$$(3.23) \qquad \frac{\partial}{\partial t} u(t,y) + \eta_t^V \frac{\partial}{\partial y} u(t,y) + \tilde{g}^y(t, u_t) = 0,$$

$$(3.24) \qquad u(T,y) = 0,$$

fulfills the equation

$$\tilde{g}^y(t, \widehat{u}_t) = g^y(t, \widehat{u}_t) \qquad \forall (t,y) \in [0,T] \times E.$$

We then conclude that $\widehat{u}$ is a weak solution to the partial differential equation (3.18) with boundary condition (3.19). In a final step, we will show that the solution is also a classical solution.

*Step 1. Definition of the auxiliary function $\tilde{g} \colon [0,T] \times \mathcal{C}_b(E) \to \mathcal{C}_b(E)$.* We introduce the function

$$\tilde{g}(t,v) := g(t, \kappa(v,t)),$$

defined on $[0,T] \times \mathcal{C}_b(E)$, with the function $\kappa$ truncating $v \in \mathcal{C}_b(E)$ in the following way. Letting $C$ be some positive constant,

$$\kappa(v,t)(x) := \max(\min(C(T-t), v(x)), -C(T-t)).$$

*Step 2. Condition* (a-2) *of Theorem* 3.6 *is fulfilled.* We have to prove that $\tilde{g}$ is a Lipschitz-continuous function on $\mathcal{C}_b(E)$, uniformly in $t$. This follows if we can show that there exists a constant $L$, independent of $(t,y) \in [0,T] \times E$, such that

$$|g^y(t, v_1) - g^y(t, v_2)| \leq L\|v_1 - v_2\|_\infty$$

for all $v_1, v_2 \in \mathcal{C}_b^Q(E)$, where $Q = CT$. In what follows, we fix a pair $(t,y) \in [0,T] \times E$ and drop the indices $(t,y)$ in the notation. We may write

$$
\begin{aligned}
|g(\cdot, v_1) - g(\cdot, v_2)| \leq \frac{1}{2(\sigma^M)^2} & \left| \left( \int W^M(x)(W^L(v_1))(x)\nu(dx) \right)^2 \right. \\
& \left. - \left( \int W^M(x)(W^L(v_2))(x)\nu(dx) \right)^2 \right| \\
& + \frac{|\int W^M(x)\nu(dx) - \widehat{\lambda} \int (W^M(x))^2 \nu(dx)|}{(\sigma^M)^2} \\
& \times \left| \int (W^L(v_1) - W^L(v_2))(x) W^M(x) \nu(dx) \right| \\
& + \left| \int (W^L(v_1) - W^L(v_2))(x) \nu(dx) \right|.
\end{aligned}
$$

By Assumption 3.1, $1/(\sigma^M)^2$ and $|\int W^M(x)\nu(dx) - \widehat{\lambda}\int (W^M(x))^2\nu(dx)|$ are uniformly bounded on $[0,T]\times E$. Moreover, $W^L(v)$ is uniformly bounded in $v\in\mathcal{C}_b^Q(E)$ by some constant $K$, and we may write, using the elementary inequality $a^2 - b^2 \le 2\max(|a|,|b|)|a-b|$,

$$
\begin{aligned}
|g(\cdot, v_1) - g(\cdot, v_2)| &\le \left(\frac{\|W^M\|_\infty}{(\sigma^M)^2}\Big[(K+1)\Big|\int W^M(x)\nu(dx)\Big| \right.\\
&\qquad\qquad\qquad + \widehat{\lambda}\int (W^M(x))^2\nu(dx)\Big] + 1\Big)\\
&\quad \times \|W^L(v_1) - W^L(v_2)\|_1.
\end{aligned}
$$

(3.25)

Due to Lemma 3.5 (Lipschitz-continuity of $W^L$), we conclude that $\tilde{g}\colon [0,T]\times \mathcal{C}_b(E) \to \mathcal{C}_b(E)$ is Lipschitz-continuous on $\mathcal{C}_b(E)$, uniformly in $t$.

Now Theorem 3.6 can be applied to the problem (3.23)–(3.24) which gives us a unique bounded weak solution $\widehat{u}\in\mathcal{C}_b([0,T]\times E)$.

*Step 3.* There exists a constant $C$ such that for all $(t,y)\in [0,T]\times E$.

$$|\widehat{u}(t,y)| \le (T-t)C. \tag{3.26}$$

Let us fix $t\in[0,T]$, $y\in E$ as well as a positive constant $C$ (to be specified below) and define [with $\hat{V}$ from (3.21)] the deterministic time $\tau_y$ as

$$\tau_y := \inf\{s\in[t,T]\ |\widehat{u}(s,\hat{V}_s^{t,y}) < (T-s)C\}\wedge T.$$

Then, $\widehat{u}(s,\hat{V}_s^{t,y}) \ge (T-s)C$, for all $s\in[t,\tau_y)$, and $\widehat{u}(\tau_y,\hat{V}_{\tau_y}^{t,y}) \le (T-\tau_y)C$. Since $\widehat{u}(s,\hat{V}_s^{t,y}) \ge (T-s)C$ for all $s\in[t,\tau_y)$, we get (with the truncation function $\kappa$ from step 1) $(\Delta\kappa(\widehat{u}_s,s))(\hat{V}_s^{t,y}) \le 0$. It follows that, for $s\in[t,\tau_y)$, the process

$$
\begin{aligned}
&W_{s,\hat{V}_s^{t,y}}^L(\kappa(\widehat{u}_s,s),x)\\
&= \exp\Big\{(\Delta\kappa(\widehat{u}_s,s))(\hat{V}_s^{t,y},x) - \Big[\widehat{\lambda} + \frac{\int W^M(z)W^L(z)\nu(dz)}{(\sigma^M)^2}\Big]W^M(x)\Big\}\\
&\quad - 1 + \widehat{\lambda}W^M(x)
\end{aligned}
$$

is bounded by some constant independent of level $C$. By our assumptions, we then can conclude from (3.22) that there exists a constant $C_1$, independent of $\tau_y$ and hence also of $C$, such that $|\tilde{g}^{\hat{V}_s^{t,y}}(s,\widehat{u}_s)| < C_1$ for all $s\in[t,\tau_y)$. It results that

$$
\begin{aligned}
\widehat{u}(t,y) &= \int_t^T \tilde{g}^{\hat{V}_s^{t,y}}(s,\widehat{u}_s)\,ds\\
&= \int_t^{\tau_y} \tilde{g}^{\hat{V}_s^{t,y}}(s,\widehat{u}_s)\,ds + \int_{\tau_y}^T \tilde{g}^{\hat{V}_s^{t,y}}(s,\widehat{u}_s)\,ds
\end{aligned}
$$

$$= \int_t^{\tau_y} \tilde{g}^{\hat{V}_s^{t,y}}(s, \hat{u}_s)\, ds + \hat{u}(\tau_y, \hat{V}_{\tau_y}^{t,y})$$

$$\leq (\tau_y - t)C_1 + (T - \tau_y)C.$$

The lower bound can be shown directly. We know that $W_s^L$ is bounded from below by $-1 + \hat{\lambda}_s W_s^M(x)$. As a direct consequence of this, together with $\sigma^M$ being bounded from above (this is the only place where we need this additional assumption), it follows that $\tilde{g}(s, \hat{u}_s)$ is bounded from below. Therefore, there exists a constant $C_2 > 0$ such that

$$\hat{u}_t \geq -(T - t)C_2.$$

If we now choose $C \geq C_1 \vee C_2$, (3.26) follows directly.

*Step 4. $\hat{u}$ is continuously differentiable in the space variable.* We use here an auxiliary function $\breve{g}$, slightly different from $\tilde{g}$. A truncation function $\breve{\kappa}$ is now introduced in such a way that we do not bound $u$, but rather the difference $\Delta u$ [i.e., we consider $\breve{\kappa}(\Delta u_t, t)$ instead of $\Delta \breve{\kappa}(u_t, t)$]. In terms of the function $\breve{g}$, this means that we work with the function

$$\breve{W}^L : \mathcal{C}_b(E) \to l^1(\operatorname{supp}(\nu))$$

defined as

$$\breve{W}^L(x) := \exp\left\{ \breve{\kappa}(\Delta u, \cdot)(x) - \left[ \hat{\lambda} + \frac{\int W^M(z) \breve{W}^L(z) \nu(dz)}{(\sigma^M)^2} \right] W^M(x) \right\}$$

$$- 1 + \hat{\lambda} W^M(x).$$

In addition, to ensure that $\breve{g}(t, u_t)$ is differentiable, we assume that $\breve{\kappa}$ has the following form, with $w \in l^\infty(\operatorname{supp}(\nu))$:

$$\breve{\kappa}(w, t)(x) = \begin{cases} v(w), & \text{if } |w(x)| \leq (T - t)C, \\ \varphi(w, t)(x), & \text{if } (T - t)C < |w(x)| < K + (T - t)C, \\ \operatorname{sign}(w(x))(K + C(T - t)), \\ & \text{if } |w(x)| \geq K + (T - t)C \end{cases}$$

for some fixed constants $C$, $K$ and a suitable $\varphi(w, t) \in l^\infty(\operatorname{supp}(\nu))$ with $|\varphi(w, t)(x)| \leq K + (T - t)C$, such that $\breve{\kappa} : l^\infty(\operatorname{supp}(\nu)) \times [0, T] \to l^\infty(\operatorname{supp}(\nu))$ is differentiable in $w$ with uniformly bounded partial derivative.

Reasoning as in Step 2, it follows that $\breve{g}$ is Lipschitz-continuous and, therefore, we may apply Theorem 3.6, which provides a solution $\breve{u}$. Let us consider $\hat{u}$ introduced above, which is bounded due to Step 3. That is, there exists a pair $(C, K)$ such that $\breve{\kappa}(\Delta \hat{u}_t, t) = \Delta \hat{u}_t$ and, therefore, $\breve{g}(t, \hat{u}_t) = g(t, \hat{u}_t)$. By uniqueness of solution, we conclude that $\breve{u} = \hat{u}$.

Let us now assume that $u_t \in \mathcal{C}_b^1(E)$. By direct calculation,

$$\frac{\partial}{\partial y} \breve{g}^y(t, u_t) = \int \left( \frac{\partial}{\partial y} \breve{\kappa}(\Delta u_{t,y})(x) \right) \tilde{W}_{t,y}^L(x) \nu(dx) + k(t, y, \breve{W}_{t,y}^L(\Delta u_{t,y}))$$

with $\tilde{W}^L(x) := \breve{W}^L(x) + 1 - \widehat{\lambda} W^M(x)$ and a uniformly bounded $k(t, y, \breve{W}^L_{t,y}(\Delta u_{t,y}))$. Let us now write

$$
\begin{aligned}
\frac{\partial}{\partial y} \breve{g}^y(t, u_t) &= \int \frac{\partial}{\partial w} \breve{\kappa}(w, t)(x) \Big|_{w = \Delta u_{t,y}} \left( \frac{\partial}{\partial y} \Delta u_{t,y}(x) \right) \tilde{W}^L_{t,y}(x) \nu(dx) \\
&\quad + k(t, y, \breve{W}^L_{t,y}(\Delta u_{t,y})) \\
&= \int \frac{\partial}{\partial w} \breve{\kappa}(v, t)(x) \Big|_{w = \int_0^{W^V_{t,y}(\cdot)} (\partial/\partial y) u(t, y+z) \, dz} \\
&\quad \times \left( \frac{\partial}{\partial y} u(t, y + W^V_{t,y}(x)) - \frac{\partial}{\partial y} u(t, y) \right) \tilde{W}^L_{t,y}(x) \nu(dx) \\
&\quad + k\left( t, y, \breve{W}^L_{t,y}\left( \int_0^{W^V_{t,y}(\cdot)} \frac{\partial}{\partial y} u(t, y+z) \, dz \right) \right) \\
&=: \hat{g}^y\left( t, \frac{\partial}{\partial y} u_t \right).
\end{aligned}
$$

Let us set $v_t(y) = \frac{\partial}{\partial y} u_t(y)$, which belongs to $\mathcal{C}_b(E)$. We already know that $\breve{W}^L$ is uniformly continuous and bounded in $(t, y, v_t) \in [0, T] \times E \times \mathcal{C}_b(E)$ and, therefore, $k$ is uniformly continuous and bounded on this set. On the other hand, taking into account the definition of $\breve{\kappa}$, it follows directly that condition (3.17) is fulfilled and that

$$
\frac{\partial}{\partial w} \breve{\kappa}(w, t)(x) \Big|_{w = \int_0^{W^V_{t,y}(\cdot)} v_t(y+z) \, dz}
$$

is uniformly continuous in $(t, y, v_t) \in [0, T] \times E \times M$. Therefore, all conditions of Theorem 3.7 are fulfilled and, hence, the solution $\widehat{u}$ to the PDE (3.18) with boundary condition (3.19) is continuously differentiable in the space variable. □

Having proved the existence of a solution to the partial differential equation (3.9) with boundary condition (3.10), we are in a position to determine the triplet $(\widehat{\phi}, W^L, \sigma^L)$ which solves equation (2.5). Since $\widehat{u}$ is uniformly bounded, we directly see that this also holds for $\widehat{\phi}$. The extra assumption that $\sigma^M$ is bounded from above is not fulfilled in some examples. We shall indicate later, using the result of Theorem 3.8, how to proceed in the standard BN–S model without this assumption and still get a uniformly bounded $\widehat{\phi}$.

THEOREM 3.9. *Let Assumption 3.1 hold, and further assume that $\sigma^M$ is uniformly bounded from above on $[0, T] \times E$. Let us assume that the triplet* $(\widehat{\phi}, W^L, \sigma^L)$ solves equation (2.5) as well as (2.3), with $(\widehat{\phi}, W^L, \sigma^L)$ uniformly bounded. Then the process $Z = (Z_t)$ defined by

$$Z_t = \frac{dQ^*}{dP}\Big|_{\mathcal{F}_t} = \mathcal{E}\left(-\left[\int (\widehat{\lambda}\sigma^M - \sigma^L) \, dY^c + (\widehat{\lambda}W^M(x) - W^L(x)) * (\mu_Y - \nu_Y)\right]\right)_t$$

is the density process of the MEMM.

PROOF.   To show that $Q^*$ is the MEMM, we show that, according to our approach outlined in Remark 2.8, $Q^*$ is an equivalent martingale measure, $I(Q^*, P) < \infty$ and $\int \frac{\widehat{\phi}}{S-} \, dS$ is a true $Q$-martingale for all $Q \in \mathcal{M}^e$ with finite relative entropy.

1. $Q^*$ is an equivalent martingale measure: Let us first show that it is an equivalent probability measure by checking the conditions of Lemma 2.11. We consider the local martingale $N$ defined by

$$
\begin{aligned}
(3.27) \quad N &= -\int \lambda \, dM + L \\
&= \int (\sigma^L - \widehat{\lambda}\sigma^M) \, dY^c + (W^L(x) - \widehat{\lambda}W^M(x)) * (\mu_Y - \nu_Y).
\end{aligned}
$$

Since $W^L$, $\widehat{\lambda}$ and $W^M$ are bounded, $N$ is locally bounded and due to

$$W^L(x) - \widehat{\lambda}W^M(x) > -1,$$

we have $\Delta N > -1$. Moreover, we set

$$
\begin{aligned}
U = \tfrac{1}{2}\int (\sigma^L - \widehat{\lambda}\sigma^M)^2 \, ds \\
+ \{(1 - \widehat{\lambda}W^M(x) + W^L(x)) \log(1 - \widehat{\lambda}W^M(x) + W^L(x)) \\
+ \widehat{\lambda}W^M(x) - W^L(x)\} * \mu_Y.
\end{aligned}
$$

Since $\widehat{\lambda}W^M$, $W^L \in l^\infty(\mathrm{supp}(\nu)) \cap l^1(\mathrm{supp}(\nu))$ and $\widehat{\lambda}$, $\sigma^M$ and $\sigma^L$ are all uniformly bounded, $U$ has locally integrable variation and its compensator $B$ is also bounded. Hence, condition (2.7) is naturally fulfilled and, therefore, $Q^*$ is an equivalent probability measure. Finally, $Q^*$ is a martingale measure since its density process can be written as

$$Z = \mathcal{E}\left(-\int \lambda \, dM + L\right),$$

where $L$ and $[M, L]$ are locally bounded local $P$-martingales.

2. $I(Q^*, P) < \infty$: The density $Z^* = \frac{dQ^*}{dP}$ may be written as

$$Z^* = \exp\left\{c + \int_0^T \frac{\widehat{\phi}_t}{S_{t-}} \, dS_t\right\},$$

where $c$ is the normalizing constant. We get

$$I(Q^*, P) = E_{Q^*}\left[c + \int_0^T \frac{\widehat{\phi}_t}{S_{t-}}\, dS_t\right]$$

$$= E_{Q^*}\left[c + \int_0^T \widehat{\phi}_t(\eta_t^M\, dt + \sigma_t^M\, dY_t^c) + (\widehat{\phi}W^M(x) * (\mu_Y - \nu_Y))_T\right].$$

We must therefore show that

$$(3.28) \qquad E_{Q^*}\left[\int_0^T \frac{\widehat{\phi}_t}{S_{t-}}\, dS_t\right] = 0,$$

since that implies $I(Q^*, P) = c$, which is finite by the previous step. Introducing

$$(3.29) \qquad \nu_Y^{Q^*} = (1 - \widehat{\lambda}W^M(x) + W^L(x)) * \nu_Y,$$

it follows from Girsanov's theorem that $W^M(x) * (\mu_Y - \nu_Y^{Q^*})$ and $\int \sigma^M\, dY^c + \int(\widehat{\lambda}\sigma^M - \sigma^L)\sigma^M\, dt$ are local $Q^*$-martingales. In fact, they are true $Q^*$-martingales since their quadratic variations are $Q^*$-integrable. This follows for the first term since $1 - \widehat{\lambda}W^M(x) + W^L(x)$ is bounded, and $W^M$ is uniformly bounded and integrable w.r.t. $\nu$. For the second term, it follows from the boundedness of $\sigma^M$. Equation (3.28) follows since the dynamics of $S$ can be written as

$$\frac{dS_t}{S_{t-}} = \sigma_t^M dY_t^c + (\widehat{\lambda}_t\sigma_t^M - \sigma_t^L)\sigma_t^M\, dt + d(W^M(x) * (\mu_Y - \nu_Y^{Q^*}))_t.$$

3. $\int \frac{\widehat{\phi}}{S_-}\, dS$ is a true $Q$-martingale for all $Q \in \mathcal{M}^e$ with finite relative entropy. In preparation for this, let us observe that for any positive constant $\alpha$ we have

$$(3.30) \qquad E[\exp\{(\alpha(W^M(x))^2 * \mu_Y)_T\}] < \infty,$$

$$(3.31) \qquad E\left[\exp\left\{\alpha \int_0^T (\sigma_t^M)^2\, dt\right\}\right] < \infty.$$

The first inequality follows from our assumption concerning $W^M$ since then

$$E[\exp\{(\alpha(W^M(x))^2 * \mu_Y)_T\}] = \exp\{((e^{\alpha(W_t^M(x))^2} - 1) * \nu_Y)_T\} < \infty$$

(see, e.g., [17], Lemma 14.39.1). Inequality (3.31) follows since $\sigma^M$ is uniformly bounded.

We have that $\int \frac{\widehat{\phi}}{S_-}\, dS$ is a local $Q$-martingale. It will be a true $Q$-martingale by Lemma 2.12 if we can show that, for some $\beta > 0$,

$$E\left[\exp\left\{\beta \int_0^T \frac{\widehat{\phi}_t^2}{S_{t-}^2}\, d[S]_t\right\}\right] < \infty.$$

We denote $k = \sup_{t \in [0,T]} \|\widehat{\phi}\|_\infty$. Let us take $\beta = \frac{1}{2k^2}$. By the Cauchy–Schwarz inequality and (3.30), (3.31) we get

$$E\left[\exp\left\{\beta \int_0^T \left(\frac{\widehat{\phi}_t}{S_{t-}}\right)^2 d[S]_t\right\}\right]$$

$$\leq E\left[\exp\left\{\frac{1}{2}\int_0^T (\sigma_t^M)^2\,dt + \frac{1}{2}((W^M(x))^2 * \mu_Y)_T\right\}\right]$$

$$< \infty. \qquad\qquad \square$$

## 4. Computing the MEMM in special cases.

4.1. *The deterministic volatility case.*  The purpose of this section is to show how we can recover some well-known results in our setup. We consider an asset process

$$\frac{dS_t}{S_{t-}} = \eta^M(t, V_t)\,dt + \sigma^M(t, V_t)\,dY_t^c + d(W^M(\,\cdot\,, V_-, x) * (\mu_Y - \nu_Y))_t,$$

$$dV_t = \eta^V(t, V_t)\,dt,$$

fulfilling Assumptions 3.1.

COROLLARY 4.1.  *Let the bounded function $\widehat{\phi}\colon [0,T] \to \mathbb{R}$ be such that*

$$(4.1) \qquad (|W^M(x)(\exp\{\widehat{\phi}W^M(x)\} - 1)| * \nu_Y)_T < \infty,$$

*and that, for any $t \in [0,T]$, the following equation is fulfilled:*

$$(4.2)\quad \eta_t^M + (\sigma_t^M)^2\widehat{\phi}_t + \int_\mathbb{R} W_t^M(x)(\exp\{\widehat{\phi}_t W_t^M(x)\} - 1)\nu_Y(dx) = 0.$$

*Then, the MEMM $Q^*$ is given by*

$$\frac{dQ^*}{dP} = \exp\left\{c + \int_0^T \frac{\widehat{\phi}_t}{S_{t-}}\,dS_t\right\}$$

*(with normalizing constant $c$). Its density process can be written as*

$$Z_t = \frac{dQ^*}{dP}\bigg|_{\mathcal{F}_t} = \mathcal{E}\left(\int \widehat{\phi}\sigma^M\,dY^c + (\exp\{\widehat{\phi}W^M(x)\} - 1) * (\mu_Y - \nu_Y)\right)_t.$$

PROOF.  In the deterministic case we have $\Delta u = 0$, since $W^V = 0$. Hence, we obtain $W^L$ immediately from (3.3) as

$$W^L(x) = \widehat{\lambda}W^M(x) - 1 + \exp\{\widehat{\phi}W^M(x)\}.$$

Equation (4.2) then follows from equation (3.11) and the definition of $\widehat{\lambda}$.  $\square$

REMARK 4.2. Equation (4.2) corresponds to a condition well known in the literature. For example, equation (3.20) in [7], condition (C) in [14], condition (4.4) in Theorem B in [11], or equation (4.30) in Theorem 4.3 of [8]. For more references containing this condition (or an equivalent form of it) we refer to [11].

### 4.2. *The orthogonal volatility case.* Let us consider the asset process

$$\frac{dS_t}{S_{t-}} = \eta^M(t, V_t)\, dt + \sigma^M(t, V_t)\, dY_t^c,$$

$$dV_t = \eta^V(t, V_t)\, dt + d(W^V(\cdot, V_-, x) * \mu_Y)_t,$$

fulfilling Assumptions 3.1 so that, in particular,

$$\widehat{\lambda} = \frac{\eta^M}{(\sigma^M)^2}$$

is bounded. Assume that $E$ is compact so that $\sigma^M$ is bounded as well. We then get the following result:

COROLLARY 4.3. *The optimal strategy is*

$$\widehat{\phi} = -\widehat{\lambda}, \tag{4.3}$$

*and the density process of the MEMM is given via*

$$W^L(t, V_{t-}, x) = \frac{v(t, V_{t-} + W_t^V(x))}{v(t, V_{t-})} - 1,$$

$$\sigma^L(t, V_{t-}) = 0,$$

*where $v$ is the classical solution of the partial differential equation*

$$\frac{\partial}{\partial t} v(t, y) + \eta_t^V \frac{\partial}{\partial y} v(t, y) - \frac{1}{2} \widehat{\lambda}_t^2 (\sigma_t^M)^2 v(t, y)$$

$$+ \int_{\mathbb{R}} (v(t, y + W_t^V(x)) - v(t, y)) \nu(dx) = 0, \tag{4.4}$$

$$v(T, y) = 1. \tag{4.5}$$

PROOF. Equation (4.3) and $\sigma^L = 0$ are direct consequences of $W^M(x) = 0$ and equation (3.11). Further, (3.12) leads to

$$W^L(t, V_{t-}, x) = \exp\{u(t, V_{t-} + W_t^V(x)) - u(t, V_{t-})\} - 1. \tag{4.6}$$

We know from Theorem 3.8 that

$$\frac{\partial}{\partial t} u(t, y) + \eta_t^V \frac{\partial}{\partial y} u(t, y) - \frac{1}{2} \widehat{\lambda}_t^2 (\sigma_t^M)^2 + \int W^L(t, y, x) \nu(dx) = 0,$$

$$u(T, y) = 0 \tag{4.7}$$

has a classical bounded solution $\hat{u}$, from which we can determine $(\hat{\phi}, W^L, \sigma^L)$ and hence the MEMM by Theorem 3.9. Using the transformation $v(t, y) = \exp u(t, y)$, we get the linear boundary problem (4.4), (4.5). □

REMARK 4.4. 1. The optimal strategy in this specific case had already been identified by Grandits and Rheinländer [16] by a conditioning argument. However, while the density of the MEMM at a fixed time $T$ has a very simple form, the corresponding density process turns out to have a more complicated structure. Becherer [2] determines the density process in a model where the volatility process switches between a finite number of states.

2. The transformation $v(t, y) = \exp u(t, y)$ is very useful here since it linearizes the partial differential equation to (4.4). However, this does not apply to the general case when the jump process directly influences the asset process. As can be seen already in the deterministic volatility case, the exponential element cannot be linearized in this way.

3. Benth and Meyer-Brandis [6] determined the MEMM for the specific case of a simplified BN–S model where no jumps occur in the price process, but with $\sigma^M$ possibly unbounded. We need the boundedness of $\sigma^M$ only if we refer to Theorem 3.8 for the existence of the IPDE (4.7). Alternatively, one could directly appeal to an existence result for the linear IPDE (4.4) and then carry out the relevant verification steps, imposing analogous conditions as in [6].

4. It follows from (3.29) that the measure $\nu_Y^Q$, where $Q$ is the MEMM, is given by $\nu_Y^Q = (W^L(x) + 1) * \nu_Y$. Since $W^L$ is specified by (4.6), in general it is a non-deterministic process and, in that case, $Y$ cannot be an additive process under $Q$. We conclude that the MEMM is not in general contained in the class of structure-preserving martingale measures as considered in [24].

4.3. *The Barndorff-Nielsen Shephard model with jumps.* In [1] the price process of a stock $S = (S_t)_{t \in [0,T]}$ is defined by the exponential $\exp\{X_t\}$ with $X = (X_t)$ satisfying

$$dX_t = (\mu + \beta \sigma_t^2) \, dt + \sigma_t \, dY_t^c + d(\rho x * \tilde{\mu}_Y)_t,$$
$$d\sigma_t^2 = -\lambda \sigma_t^2 \, dt + d(x * \tilde{\mu}_Y)_t,$$

where the parameters $\mu, \beta, \rho, \lambda$ are real constants with $\lambda > 0$ and $\rho \leq 0$, and where $\tilde{\mu}_Y$ has compensator $\tilde{\nu}_Y := \lambda \nu_Y$. In addition, $Y^d$ is assumed to be a subordinator (i.e., with positive increments only) so that we have $\text{supp}(\nu) = \mathbb{R}_+$. It can be easily shown that the process $S$ may then be written as

$$\frac{dS_t}{S_{t-}} = \left( \mu + \int (e^{\rho x} - 1) \tilde{\nu}(dx) + \sigma_t^2 (\beta + \tfrac{1}{2}) \right) dt$$
$$+ \sigma_t \, dY_t^c + d((e^{\rho x} - 1) * (\tilde{\mu}_Y - \tilde{\nu}_Y))_t.$$

The process $\sigma_t^2$ is an Ornstein–Uhlenbeck process reverting toward zero and having positive jumps given by the subordinator. An explicit representation of it is given by

$$\sigma_t^2 = \sigma_0^2 \exp\{-\lambda t\} + \int_0^t \exp\{-\lambda(t-u)\}\, dY_{\lambda u}.$$

We apply the results of Section 3 and refer for one technical step (regarding the unboundedness of $\sigma^M$) to [27]. One must pay attention to the fact that we work in this specific example with the Lévy process $Y = Y^c + \tilde{Y}^d$, where $\tilde{Y}^d = Y_\lambda^d$.

COROLLARY 4.5.   *In addition to the assumptions above, let us assume*

(4.8) $$\int_0^\infty (e^{\lambda^{-1}(\beta+1/2)^2 x} - 1)\tilde{\nu}(dx) < \infty.$$

*Let $\sigma_0^2 > 0$ be fixed and denote [noting that the integrals are well defined by (4.8) and $\rho \leq 0$]*

$$\widehat{\lambda}_t = \widehat{\lambda}_t(y) := \frac{\mu + \int(e^{\rho x} - 1)\tilde{\nu}(dx) + ye^{-\lambda t}(\beta + 1/2)}{ye^{-\lambda t} + \int(e^{\rho x} - 1)^2\tilde{\nu}(dx)}$$

*which we assume to be strictly positive. The MEMM in the case of the BN–S model is then determined as follows:*

*Let us denote*

$$g^y(t, u_t) = \tfrac{1}{2}(\sigma_t^L - \widehat{\lambda}_t e^{-1/2\lambda t}\sqrt{y})^2 + \widehat{\phi}_t\widehat{\lambda}_t e^{-\lambda t}y$$

$$+ \int[W_t^L(y,x) - (\widehat{\phi}_t + \widehat{\lambda}_t)(e^{\rho x} - 1) + \widehat{\phi}_t\widehat{\lambda}_t(e^{\rho x} - 1)^2]\tilde{\nu}(dx),$$

*where $W_t^L(y,x)$ is the solution to*

$$\exp\left\{\Delta u_t(y,x) - \left[\widehat{\lambda}_t + \frac{\int(e^{\rho z} - 1)W_t^L(y,z)\tilde{\nu}(dz)}{ye^{-\lambda t}}\right](e^{\rho x} - 1)\right\}$$

$$= 1 - \widehat{\lambda}_t(e^{\rho x} - 1) + W_t^L(y,x)$$

*and*

$$\Delta u_t(y,x) = u(t, y + e^{\lambda t}x) - u(t,y),$$

$$\widehat{\phi}_t = -\frac{\int(e^{\rho x} - 1)W_t^L(y,x)\tilde{\nu}(dx)}{ye^{-\lambda t}} - \widehat{\lambda}_t,$$

$$\sigma_t^L = -\frac{\int(e^{\rho x} - 1)W_t^L(y,x)\tilde{\nu}(dx)}{\sqrt{y}e^{-(1/2)\lambda t}}.$$

*Then, the classical solution $\hat{u}$ of the integro-PDE*

$$(4.9) \qquad \frac{\partial}{\partial t} u(t, y) + g^y(t, u_t) = 0,$$

$$(4.10) \qquad u(T, y) = 0 \qquad \forall y \in E := [\sigma_0^2, \infty)$$

*determines the MEMM via $W^L$ and $\sigma^L$:*

$$Z_t = \left. \frac{dQ^*}{dP} \right|_{\mathcal{F}_t}$$
$$= \mathcal{E}\left( \int (-\widehat{\lambda}_s \sigma_s + \sigma_s^L) \, dY_s^c + (-\widehat{\lambda}(e^{\rho x} - 1) + W^L(x)) * (\tilde{\mu}_Y - \tilde{\nu}_Y) \right)_t.$$

As $\hat{\sigma}$ is in general not bounded, we may not directly apply Theorem 3.8 to prove that there exists a classical solution $\hat{u}$ to the problem (4.9)–(4.10). Resolving this issue has turned out to be surprisingly technical and has been carried out in [27], Chapter 6.6. The existence of a solution is there constructed via an Arzela–Ascoli argument from solutions which live on compact sets (their existence is therefore guaranteed by Theorem 3.8). Let us summarize this analysis:

THEOREM 4.6 ([27]).    *Under the assumptions of Corollary 4.5, there exists a classical solution $\hat{u}$ of the integro-PDE (4.9), (4.10) such that $\Delta \hat{u}$ is bounded from above on $[0, T] \times E$.*

PROOF OF COROLLARY 4.5.    The IPDE (4.9) with boundary condition (4.10) can be derived from the results in Section 3 by making the transformation

$$\hat{\sigma}_t^2 = e^{\lambda t} \sigma_t^2$$

such that we obtain the dynamics

$$\frac{dS_t}{S_{t-}} = \left( \mu + \lambda \int (e^{\rho x} - 1) \nu(dx) + e^{-\lambda t} \hat{\sigma}_t^2 (\beta + \tfrac{1}{2}) \right) dt$$
$$+ e^{-(1/2)\lambda t} \hat{\sigma}_t \, dY_t^c + d((e^{\rho x} - 1) * (\tilde{\mu}_Y - \tilde{\nu}_Y))_t$$

with

$$d\hat{\sigma}_t^2 = d(e^{\lambda \cdot} x * \tilde{\mu}_Y)_t.$$

By Theorem 4.6, we have a classical solution $\hat{u}$, with $\Delta \hat{u}$ bounded from above, so it follows from Lemma 3.3.1 that $W^L$ is uniformly bounded and integrable w.r.t. $\nu$.

Based on this result, we now show that the three conditions of Remark 2.8 are fulfilled:

1. $Q^*$ is an equivalent martingale measure: here we proceed similarly as in the proof of Theorem 3.9 and concentrate only on verifying the conditions of Lemma 2.11. Let us consider

$$U = \tfrac{1}{2} \int \widehat{\phi}^2 \sigma^2 \, ds + W^U(x) * \tilde{\mu}_Y$$

with

$$W^U(x) := (W^L(x) + 1 - \widehat{\lambda}(e^{\rho x} - 1)) \log(W^L(x) + 1 - \widehat{\lambda}(e^{\rho x} - 1))$$
$$+ \widehat{\lambda}(e^{\rho x} - 1) - W^L(x).$$

Since $e^{\rho x} - 1$ and $W^L$ are uniformly bounded and integrable w.r.t. $\nu$, $U$ has locally integrable variation, and we get

$$E[\exp\{2\lambda(W^U(x) * \tilde{\nu})_T\}] < \infty.$$

Hence (2.7) is fulfilled by the Cauchy–Schwarz inequality if we can show that

$$E\left[\exp\left\{\int_0^T \widehat{\phi}_t^2 \sigma_t^2 \, dt\right\}\right] < \infty.$$

By definition, we have

$$\widehat{\phi}_t = -\widehat{\lambda}_t - \frac{\int (e^{\rho x} - 1) W_t^L(x) \tilde{\nu}(dx)}{\sigma_{t-}^2}.$$

Since $\widehat{\lambda}$ is positive and $W^L$ is bounded, $\widehat{\phi}_t$ is negative for sufficiently large $\sigma_{t-}$. Let us introduce $\overline{\sigma}$ such that for all $t \in [0, T]$ (with the possible exception of a Lebesgue-zero set)

$$\widehat{\phi}_t < 0 \qquad \text{for all } \sigma_t > \overline{\sigma}.$$

On the other hand, since $W_t^L(x) \geq -1 + \widehat{\lambda}_t(e^{\rho x} - 1)$, $\widehat{\phi}_t$ is bounded from below with

$$\widehat{\phi}_t \geq -\widehat{\lambda}_t - \frac{\int (e^{\rho x} - 1)(-1 + \widehat{\lambda}_t(e^{\rho x} - 1)) \tilde{\nu}(dx)}{\sigma_{t-}^2}$$

$$= -\left(\beta + \frac{1}{2}\right) - \frac{\mu}{\sigma_{t-}^2}$$

because of

$$\widehat{\lambda}_t = \frac{\mu + \int (e^{\rho x} - 1) \tilde{\nu}(dx) + \sigma_{t-}^2 (\beta + 1/2)}{\sigma_{t-}^2 + \int (e^{\rho x} - 1)^2 \tilde{\nu}(dx)}.$$

Let us now analyze

$$E\left[\exp\left\{\int_0^T \widehat{\phi}_t^2 \sigma_t^2 \, dt\right\}\right]$$

$$= E\left[\exp\left\{\int_0^T \mathbf{1}_{\{\sigma_t \leq \overline{\sigma}\}} \widehat{\phi}_t^2 \sigma_t^2 \, dt\right\} \exp\left\{\int_0^T \mathbf{1}_{\{\sigma_t > \overline{\sigma}\}} \widehat{\phi}_t^2 \sigma_t^2 \, dt\right\}\right].$$

We have that

$$\exp\left\{\int_0^T \mathbf{1}_{\{\sigma_t \le \overline{\sigma}\}} \widehat{\phi}_t^2 \sigma_t^2 \, dt\right\}$$

is uniformly bounded. Moreover,

$$E\left[\exp\left\{\int_0^T \mathbf{1}_{\{\sigma_t > \overline{\sigma}\}} \widehat{\phi}_t^2 \sigma_t^2 \, dt\right\}\right]$$

is finite due to (i) the fact that for almost all $t$

$$0 \ge \widehat{\phi}_t \ge -\left(\beta + \frac{1}{2}\right) - \frac{\mu}{\sigma_t^2}$$

on the set $\{\sigma_t > \overline{\sigma}\}$ and (ii) condition (4.8), which, according to Benth, Karlsen and Reikvam ([5], Lemma 3.1), ensures that

$$E\left[\exp\left\{(\beta + \tfrac{1}{2})^2 \int_0^T \sigma_t^2 \, dt\right\}\right] < \infty.$$

2. $I(Q^*, P) < \infty$: We have to show that for $\tilde{\nu}_Y^{Q^*} = (W^L(x) + 1 - \widehat{\lambda}(e^{\rho x} - 1)) * \tilde{\nu}_Y$,

$$(e^{\rho x} - 1) * (\tilde{\mu}_Y - \tilde{\nu}_Y^{Q^*}) \quad \text{and} \quad \int \sigma \, dY^c + \int (\widehat{\lambda}\sigma - \sigma^L)\sigma \, dt$$

are true $Q^*$-martingales, that is, their quadratic variations are $Q^*$-integrable. This follows for the first term from the boundedness of $W^L$ and the integrability of $e^{\rho x} - 1$. For the second term, let us consider

$$E_{Q^*}\left[\left[\int \sigma \, dY^c\right]_T\right] = E_{Q^*}\left[\int_0^T \sigma_t^2 \, dt\right].$$

It is well known that we may write

$$\int_0^T \sigma_t^2 \, dt = \lambda^{-1}(1 - e^{-\lambda T})\sigma_0^2 + (\lambda^{-1}(1 - e^{-\lambda(T - \cdot)})x * \tilde{\mu}_Y)_T.$$

Hence, $E_{Q^*}[\int_0^T \sigma_t^2 \, dt]$ is finite if $E_{Q^*}[(x * \tilde{\nu}_Y^{Q^*})_T]$ is finite, which, since $W^L$ is bounded, is equivalent to showing that $\int x\nu(dx) < \infty$. However, this follows from condition (4.8).

3. $\int \frac{\widehat{\phi}}{S_-} \, dS$ is a true $Q$-martingale for all $Q \in \mathcal{M}^e$ with finite relative entropy: by Lemma 2.12, $\int \frac{\widehat{\phi}}{S_-} \, dS$ is a true $Q$-martingale if we can show that, for some $\gamma > 0$,

$$E\left[\exp\left\{\gamma \int \frac{\widehat{\phi}_t^2}{S_{t-}^2} \, d[S]_t\right\}\right]$$

$$= E\left[\exp\left\{\gamma \int_0^T \widehat{\phi}_t^2 \sigma_t^2 \, dt + (\gamma \widehat{\phi}^2(e^{\rho x} - 1)^2 * \tilde{\mu}_Y)_T\right\}\right]$$

$$< \infty.$$

We have

$$E[\exp\{(2\gamma\widehat{\phi}^2(e^{\rho x}-1)^2 * \tilde{\mu}_Y)_T\}] < \infty,$$

and, for $\gamma < \frac{\beta+1/2}{2\max\tilde{\phi}_t^2}$, it follows that

$$E\left[\exp\left\{2\gamma\int_0^T \widehat{\phi}_t^2\sigma_t^2\,dt\right\}\right] < \infty.$$

Therefore, an application of the Cauchy–Schwarz inequality yields

$$E\left[\exp\left\{\gamma\int_0^T \frac{\widehat{\phi}_t^2}{S_{t-}^2}\,d[S]_t\right\}\right] < \infty. \qquad\qquad\square$$

## APPENDIX

PROOF OF LEMMA 3.3. We consider the equation

$$\Phi = \int f(x)\exp\{k(x)-\beta f(x)\Phi\}\nu(dx)$$

and will show that there exists a unique value $\Phi_k \in \mathbb{R}$ which solves it. For this purpose, let us define

$$H(z) = z - \int f(x)\exp\{k(x)-\beta f(x)z\}\,\nu(dx).$$

Since

$$\lim_{z\to\infty} -f(x)\exp\{-\beta f(x)z\} = \begin{cases} 0, & f(x) \geq 0, \\ \infty, & f(x) < 0, \end{cases}$$

we have $\lim_{z\to\infty} H(z) = \infty$ and, for reasons of symmetry, $\lim_{z\to-\infty} H(z) = -\infty$. Furthermore, $H$ is continuously differentiable with

$$\frac{\partial}{\partial z}H(z) = 1 + \int \beta f^2(x)\exp\{k(x)-\beta f(x)z\}\nu(dx) > 0.$$

Therefore, there exists a unique $\Phi_k \in \mathbb{R}$ such that $H(\Phi_k) = 0$. We can moreover show that

$$|\Phi_k| \leq \max_{x\in\mathrm{supp}(\nu)}\{\exp k(x)\}\int |f(x)|\nu(dx).$$

Let us assume that $\Phi_k \geq 0$. Then, we get

$$\Phi_k = \int f(x)\exp\{k(x)-\beta f(x)\Phi_k\}\nu(dx)$$

$$\leq \int_{\{f(x)>0\}} f(x)\exp\{k(x)-\beta f(x)\Phi_k\}\nu(dx)$$

$$\leq \int_{\{f(x)>0\}} f(x) \exp\{k(x)\} \nu(dx)$$

$$\leq \max_{x\in\mathrm{supp}(\nu)} \{\exp k(x)\} \int_{\{f(x)>0\}} f(x)\nu(dx)$$

$$\leq \max_{x\in\mathrm{supp}(\nu)} \{\exp k(x)\} \int |f(x)|\nu(dx).$$

The lower bound can be shown in exactly the same way.

Let us now define the bounded function

$$\varphi_k(x) := \exp\{k(x) - \beta f(x)\Phi_k\}.$$

As we have

$$\int f(x)\varphi_k(x)\nu(dx) = \int f(x)\exp\{k(x) - \beta f(x)\Phi_k\}\nu(dx) = \Phi_k,$$

it follows that

$$\varphi_k(x) = \exp\left\{k(x) - \beta f(x)\int f(z)\varphi_k(z)\nu(dz)\right\},$$

and, therefore, we conclude that $\varphi := \varphi_k$ is well defined and bounded.  $\square$

PROOF OF LEMMA 3.5.  Since $W^L$ is bounded on $\mathcal{C}_b^Q(E)$, we only have to show local Lipschitz-continuity of $W^L$, that is, we have to show that for any $c > 0$, there exists a constant $L_c$ such that

$$\|W^L(v_1) - W^L(v_2)\|_1 \leq L_c\|v_1 - v_2\|_\infty$$

for all $v_1, v_2 \in \mathcal{C}_b^Q(E)$ with $\|v_1 - v_2\|_\infty \leq \frac{c}{2}$. For that purpose, consider $v_0 + rh$, where $v_0 \in \mathcal{C}_b^Q(E)$ and $h \in \mathcal{C}_b^Q(E)$ with $\|h\|_\infty = \frac{c}{2}$, and $r \in [0,1]$.

Let $k$ be bounded from above and define

$$\varphi_k(x) := \exp\left\{k(x) - \frac{W^M(x)}{(\sigma^M)^2}\int W^M(z)\varphi_k(z)\nu(dz)\right\},$$

$$\Phi_k := \int W^M(x)\varphi_k(x)\nu(dx).$$

By equation (3.12), we may write

$$\varphi_{k^*(r)} = W^L(v_0 + rh) - \widehat{\lambda}W^M + 1$$

for

$$k^*(r) = \Delta v_0 + r\Delta h - W^M\widehat{\eta}$$

and

$$\widehat{\eta} := \widehat{\lambda}\left(1 + \frac{\int (W^M(z))^2\nu(dz)}{(\sigma^M)^2}\right) - \frac{\int W^M(z)\nu(dz)}{(\sigma^M)^2}.$$

The goal is to show that there is a constant $C_1$ such that, for all $r \in [0,1]$,

$$(A.1) \qquad \|\varphi_{k^*(r)} - \varphi_{k^*(0)}\|_1 \leq C_1 \|rh\|_\infty = \frac{C_1 rc}{2}.$$

Let us therefore analyze

$$
\begin{aligned}
& |(\varphi_{k^*(r)} - \varphi_{k^*(0)})(x)| \\
& \quad = \exp\{\Delta v_0(x) - \widehat{\eta} W^M(x)\} \\
(A.2) & \qquad \times \left| \exp\left\{ r\Delta h(x) - \Phi_{k^*(r)} \frac{W^M(x)}{(\sigma^M)^2} \right\} - \exp\left\{ -\Phi_{k^*(0)} \frac{W^M(x)}{(\sigma^M)^2} \right\} \right| \\
& \quad = \exp\left\{ \Delta v_0(x) - \widehat{\eta} W^M(x) - \Phi_{k^*(0)} \frac{W^M(x)}{(\sigma^M)^2} \right\} \\
& \qquad \times \left| \exp\left\{ r\Delta h(x) - (\Phi_{k^*(r)} - \Phi_{k^*(0)}) \frac{W^M(x)}{(\sigma^M)^2} \right\} - 1 \right|.
\end{aligned}
$$

Since $v_0$ is uniformly bounded by $Q$, the first term on the right-hand side is uniformly bounded for all $x \in \operatorname{supp}(\nu)$. The second term [to be labeled $f_x(r)$] needs further investigation. For this purpose, let us state the following property of $\Phi_k$:

CLAIM 1. *Given two functions $k_1, k_2 \in l^\infty(\operatorname{supp}(\nu))$ with*

$$
\begin{cases}
k_1(x) \leq k_2(x) & \forall x \in \operatorname{supp}(\nu) \ s.t. \ W^M(x) < 0, \\
k_1(x) \geq k_2(x) & \forall x \in \operatorname{supp}(\nu) \ s.t. \ W^M(x) > 0,
\end{cases}
$$

*it follows that $\Phi_{k_1} \geq \Phi_{k_2}$.*

PROOF. Let us assume that $\Phi_{k_1} < \Phi_{k_2}$. For any $x \in \operatorname{supp}(\nu)$, it then follows that

$$\frac{\varphi_{k_2}(x)}{\varphi_{k_1}(x)} = \exp\left\{ k_2(x) - k_1(x) - \frac{W^M(x)}{(\sigma^M)^2}(\Phi_{k_2} - \Phi_{k_1}) \right\}$$

$$
\begin{cases}
> 1, & \forall x \in \operatorname{supp}(\nu) \ \text{s.t.} \ W^M(x) < 0, \\
< 1, & \forall x \in \operatorname{supp}(\nu) \ \text{s.t.} \ W^M(x) > 0.
\end{cases}
$$

However, this leads to a contradiction, since then

$$\Phi_{k_2} - \Phi_{k_1} = \int W^M(x)(\varphi_{k_2}(x) - \varphi_{k_1}(x))\nu(dx) \geq 0.$$

Therefore, we must have $\Phi_{k_1} \geq \Phi_{k_2}$. $\quad \square$

Let us now fix $x_0 \in \operatorname{supp}(\nu)$ and analyze the term

$$f_{x_0}(r) := \exp\left\{ r\Delta h(x_0) - (\Phi_{k^*(r)} - \Phi_{k^*(0)}) \frac{W^M(x_0)}{(\sigma^M)^2} \right\} - 1.$$

Obviously, we have $f_{x_0}(0) = 0$. Let us now assess the upper and lower bounds of $f_{x_0}$ for $r \in [0,1]$. For this purpose, we introduce

$$k^+(r,x) := \Delta v_0(x) - W^M(x)\widehat{\eta} + rc(\mathbf{1}_{\{W^M(x)>0\}} - \mathbf{1}_{\{W^M(x)<0\}})(x),$$

$$k^-(r,x) := \Delta v_0(x) - W^M(x)\widehat{\eta} - rc(\mathbf{1}_{\{W^M(x)>0\}} - \mathbf{1}_{\{W^M(x)<0\}})(x).$$

We will use in the following the notation

$$\mathbf{1}^*(x) := \mathbf{1}_{\{W^M(x)>0\}} - \mathbf{1}_{\{W^M(x)<0\}}.$$

It follows from the claim above that

(A.3)          $\Phi_{k^-(r)} - \Phi_{k^-(0)} \leq \Phi_{k^*(r)} - \Phi_{k^*(0)} \leq \Phi_{k^+(r)} - \Phi_{k^+(0)}.$

Let us now consider in detail the upper bound,

$$\Phi_{k^+(r)} - \Phi_{k^+(0)} = \int_0^r \frac{\partial}{\partial r} \Phi_{k^+(s)} \, ds.$$

Here the existence of the derivative can be guaranteed by an application of the Implicit Function Theorem for Banach spaces (see, e.g., [28], page 150) to the equation

$$\Phi_{k^+(r)} = \int W^M(x) \exp\left\{ k^+(r,x) - \frac{W^M(x)}{(\sigma^M)^2}\Phi_{k^+(r)} \right\} \nu(dx).$$

We have

$$\frac{\partial}{\partial r}\Phi_{k^+(r)} = \int W^M(x)\left[ c\mathbf{1}^*(x) - \frac{W^M(x)}{(\sigma^M)^2}\frac{\partial}{\partial r}\Phi_{k^+(r)} \right]$$

$$\times \exp\left\{ k^+(r,x) - \frac{W^M(x)}{(\sigma^M)^2}\Phi_{k^+(r)} \right\} \nu(dx),$$

so we may write [recalling that $\varphi_{k^+(r)}(x) = \exp\{k^+(r,x) - \frac{W^M(x)}{(\sigma^M)^2}\Phi_{k^+(r)}\}$]

$$\frac{\partial}{\partial r}\Phi_{k^+(r)} = c\left( \frac{\int W^M(x)\mathbf{1}^*(x)\varphi_{k^+(r)}(x)\nu(dx)}{1 + \int ((W^M(x))^2/(\sigma^M)^2)\varphi_{k^+(r)}(x)\nu(dx)} \right)$$

$$< c\int |W^M(x)|\varphi_{k^+(r)}(x)\nu(dx).$$

Since $k^+(s)$ is bounded from above, it follows from the definition of $\varphi_{k^+(s)}$ and Lemma 3.3 that $\varphi_{k^+(s)}$ is uniformly bounded by some constant $K^*$ for any $s \in [0,r]$. Therefore, it follows that

$$\Phi_{k^+(r)} - \Phi_{k^+(0)} < crK^* \int |W^M(x)|\nu(dx).$$

Applying the same steps to the lower bound, it follows that

$$\Phi_{k^-(r)} - \Phi_{k^-(0)} > -crK^* \int |W^M(x)|\nu(dx).$$

Taking into account the inequalities of (A.3), we obtain the following bounds:

$$\exp\{rc\tilde{K}(x_0)\} - 1 \geq f_{x_0}(r) \geq \exp\{-rc\tilde{K}(x_0)\} - 1$$

with

$$\tilde{K}(x_0) := 1 + K^* \frac{|W^M(x_0)|}{(\sigma^M)^2} \int |W^M(x)| \nu(dx).$$

Therefore, for $r \in [0, 1]$, it follows that

$$|f_x(r)| \leq rc\tilde{K}(x)$$

with $\tilde{K} \in l^1(\text{supp}(\nu))$, and hence, via (A.1), the Lipschitz-continuity of $W^L$ is shown.  □

PROOF OF THEOREM 3.6.  Let us fix some $u \in \mathcal{C}_b([0, T] \times E)$ and consider the PDE

(A.4)  $$\frac{\partial}{\partial t} w(t, z) + b(t, z) \frac{\partial}{\partial z} w(t, z) + g^z(t, u_t) = 0,$$

(A.5)  $$w(T, z) = h(z) \qquad \forall z \in E.$$

It is straightforward to see that

$$w(t, z) = h(Z_T^{t,z}) + \int_t^T g^{Z_s^{t,z}}(s, u_s)\, ds$$

solves the boundary problem in the weak sense. Let us introduce the operator $F: \mathcal{C}_b([0, T] \times E) \to \mathcal{C}_b([0, T] \times E)$ defined as follows:

$$(Fu)(t, z) = h(Z_T^{t,z}) + \int_t^T g^{Z_s^{t,z}}(s, u_s)\, ds.$$

We have to prove that $F$ is a contraction on the space $\mathcal{C}_b([0, T] \times E)$. Let us, for some $\beta \in \mathbb{R}_+$, consider the norm

$$\|u\|_\beta := \sup_{(t,z) \in [0,T] \times E} e^{-\beta(T-t)} |u(t, z)|,$$

which is equivalent to the supremum-norm $\|u\|_\infty$. Due to condition (a-2), we obtain for $u_1, u_2 \in \mathcal{C}_b([0, T] \times E)$

$$e^{-\beta(T-t)} |(Fu_1)(t, z) - (Fu_2)(t, z)|$$

$$= \frac{1}{e^{\beta(T-t)}} \left| \int_t^T \left( g^{Z_s^{t,z}}(s, u_{1,s}) - g^{Z_s^{t,z}}(s, u_{2,s}) \right) ds \right|$$

$$\leq \frac{1}{e^{\beta(T-t)}} \int_t^T |g^{Z_s^{t,z}}(s, u_{1,s}) - g^{Z_s^{t,z}}(s, u_{2,s})| e^{-\beta(T-s)} e^{\beta(T-s)}\, ds$$

$$\leq \frac{1}{e^{\beta(T-t)}} L \|u_1 - u_2\|_\beta \int_t^T e^{\beta(T-s)}\, ds$$

$$\leq \frac{L}{\beta} \|u_1 - u_2\|_\beta$$

for all $t \in [0, T]$ and $z \in E$. Thus,

$$\|(Fu_1)(t, z) - (Fu_2)(t, z)\|_\beta \le \frac{L}{\beta} \|u_1 - u_2\|_\beta,$$

and $F$ is a contraction on the normed space $(\mathcal{C}_b([0, T] \times E), \|\cdot\|_\beta)$ with $\beta > L$. Therefore, there exists a unique fixed point $u \in \mathcal{C}_b([0, T] \times E)$ which satisfies the PDE (3.15)–(3.16) in the weak sense. $\square$

PROOF OF THEOREM 3.7.   Let us analyze the operator $\hat{G} : \mathcal{C}_b([0, T] \times E) \to \mathcal{C}_b([0, T] \times E)$, defined by

$$(\hat{G}v)(t, z) = \frac{\partial}{\partial z} h(Z_T^{t,z}) + \int_t^T \left( \frac{\partial}{\partial z} Z_s^{t,z} \right) \hat{g}^{Z_s^{t,z}}(s, v_s) \, ds.$$

Let us first discuss $\frac{\partial}{\partial z} Z_s^{t,z}$, which is well defined by Protter [25], Theorem V.39. Differentiating (3.14), we get

$$\frac{\partial}{\partial z} Z_s^{t,z} = 1 + \int_t^s \left( \frac{\partial}{\partial z} Z_u^{t,z} \right) \frac{\partial}{\partial Z_u^{t,z}} b(u, Z_u^{t,z}) \, du.$$

By Gronwall's lemma, we can directly conclude that $\frac{\partial}{\partial z} Z_s^{t,z}$ is uniformly bounded, the bound being denoted by $L_Z$. Analogously, let us denote $L_h := \|h'\|_\infty$.

Let us now discuss, for $v \in \mathcal{C}_b([0, T] \times E)$,

$$e^{-\beta(T-t)} |(\hat{G}v)(t, z)|$$
$$\le e^{-\beta(T-t)} \left( \left| \frac{\partial}{\partial z} Z_s^{t,z} \right| |h'(Z_T^{t,z})| + \int_t^T \left| \frac{\partial}{\partial z} Z_s^{t,z} \right| |\hat{g}^{Z_s^{t,z}}(s, v_s)| \, ds \right)$$
$$\le e^{-\beta(T-t)} L_Z \left( L_h + \int_t^T (L \|v_s\|_\infty + K) e^{-\beta(T-t)} e^{\beta(T-t)} \, ds \right)$$
$$\le \frac{L_Z L}{\beta} \|v\|_\beta + L_Z K T + L_Z L_h.$$

Hence, for $\beta = 2L_Z L$ and $N := \{v \in \mathcal{C}_b([0, T] \times E) | \ \|v\|_\beta \le 2L_Z(KT + L_h)\}$, $\hat{G}$ maps $N$ into $N$. Using the Arzela–Ascoli theorem, one can show that $\hat{G}$ is a compact operator on $N$. By Schauder's Fixed Point Theorem, we conclude that $\hat{G} : N \to N$ has at least one fixed point $\hat{v}$. Let us assume that $u \in \mathcal{C}_b([0, T] \times E)$ is differentiable in the space variable. Hence,

$$\frac{\partial}{\partial z} (Fu)(t, z) = \frac{\partial}{\partial z} h(Z_T^{t,z}) + \int_t^T \left( \frac{\partial}{\partial z} Z_s^{t,z} \right) \frac{\partial}{\partial Z_s^{t,z}} g^{Z_s^{t,z}}(s, u_s) \, ds$$
$$= \left( \hat{G} \frac{\partial}{\partial z} u \right)(t, z).$$

Let us now consider the primitive with respect to $z \in E$ of $\widehat{v}$, denoted $\widehat{u}$. We may write

$$\frac{\partial}{\partial z}(F\widehat{u})(t,z) = (\widehat{G}\widehat{v})(t,z) = \widehat{v}(t,z) = \frac{\partial}{\partial z}\widehat{u}(t,z).$$

It therefore follows that the function $\widehat{u}$ may be written as

$$\widehat{u}(t,z) = (F\widehat{u})(t,z) + C(t)$$

with function $C\colon [0,T] \to \mathbb{R}$. On the other hand, we know by Theorem 3.6 that there exists a unique fixed point of operator $F$ in $\mathcal{C}_b([0,T] \times E)$. Hence, choosing $C \equiv 0$, it follows that $\widehat{u}$ is uniquely defined. We have therefore shown that there exists a unique classical solution to the boundary problem (3.15)–(3.16). $\quad\square$

DEPARTMENT OF STATISTICS
LONDON SCHOOL OF ECONOMICS
LONDON WC2A 2AE
UNITED KINGDOM
E-MAIL: t.rheinlander@lse.ac.uk

D-MATH, ETH-ZENTRUM
ETH ZÜRICH
CH-8092 ZÜRICH
SWITZERLAND
E-MAIL: gallus.steiger@math.ethz.ch